\documentclass{article}
\usepackage[all]{xy}
\usepackage{amssymb}
\usepackage{amsmath}
\usepackage{amsthm}
\usepackage{tensor}
\usepackage{graphicx}
\usepackage{rotating}
\usepackage{bbm}
\usepackage{placeins}
\usepackage[T1]{fontenc}
\usepackage{pgf}
\usepackage{tikz}
\usetikzlibrary{matrix,arrows}

\DeclareMathOperator{\Hom}{Hom}

\DeclareMathOperator{\inj}{inj}

\DeclareMathOperator{\rep}{rep}
\DeclareMathOperator{\GL}{GL}
\newcommand{\G}{G}

\DeclareMathOperator{\Iso}{Iso}

\newcommand{\SN}{S}
\newcommand{\la}{\langle}
\newcommand{\ra}{\rangle}
\newcommand{\rar}{\rightarrow}

\newcommand{\xrar}{\xrightarrow}

\newcommand{\Orb}{\mathcal{O}}

\newcommand{\B}{\mathcal{B}}
\newcommand{\Pa}{\mathcal{P}}
\newcommand{\N}{\mathcal{N}}

\newcommand{\Q}{\mathcal{Q}}
\newcommand{\Qtil}{\widetilde{\Q}}
\newcommand{\dimv}{\underline{\dim}}
\newcommand{\df}{{\underline{d}}}
\newcommand{\vect}{(1,2,\ldots,n)}
\newcommand{\reps}{\rep_k(\Q,I)}
\newcommand{\repinj}{\rep_k^{\inj}(\Q,I)}

\newcommand{\U}{\mathcal{U}}
\newcommand{\V}{\mathcal{V}}
\newcommand{\W}{\mathcal{W}}

\newcommand{\remark}{{\bf Remark:  }}

\newcommand{\example}{{\bf Example:  }}

\newtheorem{theorem}{Theorem}[section]
\newtheorem{lemma}[theorem]{Lemma}

\newtheorem{proposition}[theorem]{Proposition}
\newtheorem{corollary}[theorem]{Corollary}

\begin{document}

\parindent0pt
\title{\bf $B$-orbits of $2$-nilpotent matrices and generalizations}

\author{Magdalena Boos and Markus Reineke\\ Fachbereich C - Mathematik\\ Bergische Universit\"at Wuppertal\\ D - 42097 Wuppertal\\ boos@math.uni-wuppertal.de, reineke@math.uni-wuppertal.de}
\date{}
\maketitle

\begin{abstract} The orbits of the group $B_n$ of upper-triangular matrices acting on $2$-nilpotent complex matrices via conjugation are classified via oriented link patterns, generalizing A. Melnikov's classification of the $B_n$-orbits on upper-triangular such matrices. The orbit closures as well as the ``building blocks'' of minimal degenerations of orbits are described. The classification uses the theory of representations of finite-dimensional algebras. Furthermore, we initiate the study of the $B_n$-orbits on arbitrary nilpotent matrices.
\end{abstract}

\section{Introduction}\label{intro}
The study of adjoint actions and variants thereof, and in particular the classification of orbits for such actions and the description of the orbit closures, are a common theme in Lie representation theory. The archetypical example is the Jordan-Gerstenhaber theory for the conjugacy classes of complex $n\times n$-matrices.\\[1ex]
A more recent case is A. Melnikov's study of the action of the Borel subgroup $B_n$ acting on upper-triangular $2$-nilpotent matrices via conjugation \cite{Melnikov1,Melnikov2}. The orbits and their closures are described there combinatorially in terms of so-called link patterns, which we will recapitulate in section \ref{rom}.\\[1ex]
Our aim in this paper is to generalize the work of A. Melnikov by extending the variety of upper-triangular $2$-nilpotent matrices to all $2$-nilpotent matrices. The basic setup to reach this goal is a translation of the classification problem to a problem in  representation theory of finite-dimensional algebras. More precisely, this translation yields a bijection between the orbits and the isomorphism classes of certain representations of a specific finite-dimensional algebra, see section \ref{trans}. After a brief summary of methods from the representation theory of algebras (see, for example, \cite{ASS}) in section \ref{roa}, we are able to calculate all indecomposable representations using Auslander-Reiten theory \cite{AR} in section \ref{covsection} and to classify the required representations. This gives a combinatorial classification in terms of oriented link patterns in section \ref{cob}.\\[1ex]
Since several results on orbit closures for representations of finite-dimensional algebras are available through work of G. Zwara \cite{Zwara1,Zwara2}, we can also characterize the orbit closures of $2$-nilpotent matrices in section \ref{close}.\\[1ex]
Finally, we study the conjugation action of upper-triangular matrices on arbitrary nilpotent matrices. We provide a generic normal form for the orbits of this action in section \ref{gnf} and construct a large class of semiinvariants in section \ref{si}.\\[3ex]
{\bf Acknowledgments:} The authors would like to thank K. Bongartz and A. Melnikov for valuable discussions concerning the methods and results of this work.

\section{The basic setup}\label{setup}

In this section, we fix some notation and collect information about the aforementioned group action. In addition, we summarize  material from the representation theory of finite-dimensional algebras.\\[1ex]
Let $k=\bf C$ be the field of complex numbers. We denote by $B_n\subset\GL_n(k)$ the Borel subgroup of upper-triangular matrices, by $\N_n\subset M_{n\times n}(k)$ the variety of nilpotent $n\times n$-matrices, and by $\N_n^{(2)}$ the closed subvariety of $2$-nilpotent such matrices. Obviously, $\GL_n(k)$ and $B_n$ act on $\N_n$ and on $\N_n^{(2)}$ via conjugation.\\[1ex]
In case of the action of $\GL_n(k)$ on $\N_n$, the classical Jordan-Gerstenhaber theory gives a complete classification of the orbits and their closures in terms of partitions (or, equivalently, Young diagrams).\\[1ex]
Our aim is to classify the orbits $\Orb_A$ of $2$-nilpotent matrices $A\in\N_n^{(2)}$ under the action of $B_n$. Such a classification will be given in terms of {\it oriented link patterns}; these are oriented graphs on the set of vertices $\{1,\ldots,n\}$ such that every vertex is incident with at most one arrow. This is followed by a description of the orbit closures, by giving a necessary and sufficient condition to decide whether one orbit is contained in the closure of another, and by a method to construct all orbits contained in a given orbit closure. These descriptions are also given in terms of oriented link patterns.

\subsection{Results of A. Melnikov}\label{rom}

The group $B_n$ also acts on $\mathfrak{n}_n\subset\N_n$, the space of all upper-triangular matrices in $\N_n$, and on $\mathfrak{n}_n^{(2)}=\mathfrak{n}_n\cap\N_n^{(2)}$. The orbits and their closures for the latter action are described by A. Melnikov in \cite{Melnikov1,Melnikov2}. Since these results will be generalized in the following, we describe them in more detail.\\
Let $\SN_n^{(2)}$ be the set of involutions in the symmetric group $\SN_n$ in $n$ letters. An element $\sigma$ of $S_n^{(2)}$ is represented by a so-called link pattern, an unoriented graph with vertices $\{1,\ldots,n\}$ and an edge between $i$ and $j$ if $\sigma(i)=j$.
 For example, the involution $(1,2)(3,5)\in\SN_5(k)$ corresponds to the link pattern
\begin{center}
\begin{picture}(100,30)
\multiput(10,10)(20,0){5}%
{\circle*{2}}
 \put(10,-5){1}
 \put(30,-5){2}
 \put(50,-5){3}
 \put(70,-5){4}
 \put(90,-5){5}
 \qbezier(10,10)(20,40)(30,10)
 \qbezier(50,10)(70,40)(90,10)
\end{picture}.
\end{center}
For $\sigma\in\SN_n^{(2)}$, define $N_\sigma\in\N_n^{(2)}$ by
$$(N_\sigma)_{i,j}=\left\{\begin{array}{cc}1&\mbox{ if $i<j$ and $\sigma(i)=j$,}\\ 0&\mbox{otherwise,}\end{array}\right.$$
and denote by $\B_\sigma=B_n\cdot N_\sigma$ the $B_n$-orbit of $N_\sigma$.
\begin{theorem}\label{mel1}\cite{Melnikov1} Every orbit of $B_n$ in $\mathfrak{n}_n^{(2)}$ is of the form $\B_\sigma$ for a unique $\sigma\in\SN_n^{(2)}$.
\end{theorem}
The next step is to look at the (Zariski-)closures of the orbits $\B_\sigma$.\\
For $1\leq i<j\leq n$, consider the canonical projection
$\pi_{i,j}:\mathfrak{n}_n^{(2)}\rightarrow\mathfrak{n}_{j-i+1}^{(2)}$ deleting the first $i-1$ and the last $n-j$ columns
and rows of a matrix in $\mathfrak{n}_n^{(2)}$. Define the rank matrix $R_u$ of $u\in\mathfrak{n}_n^{(2)}$ by
$$ (R_u)_{i,j}=\left\{ \begin{array}{ll}
{\rm rank}\,(\pi_{i,j}(u))&{\rm if}\ i< j;\\
                      0 &{\rm otherwise}.\\
                       \end{array}\right. $$
The rank matrix $R_u$ is $B_n$-invariant, and we denote $R_\sigma=R_{N_\sigma}$ for $\sigma\in\SN_n^{(2)}$. We define a partial ordering on the set of rank matrices by $R_{\sigma'}\preccurlyeq R_{\sigma}$ if $(R_{\sigma'})_{i,j}\leq (R_{\sigma})_{i,j}$ for all $i$ and $j$, inducing a partial ordering on $\SN_n^{(2)}$ by $\sigma'\preccurlyeq\sigma$ if $R_{\sigma'}\preccurlyeq R_{\sigma}$.
\begin{theorem}\cite{Melnikov2} The orbit closure of $\B_{\sigma}$ is given by $\overline{\B_{\sigma}}=\bigcup_{\sigma'\preccurlyeq\sigma}\B_{\sigma'}$. Moreover, the entry $(R_{\sigma})_{i,j}$ of the rank matrix
equals the number of edges with end points $e_1$ and $e_2$ such that $i\leq e_1,e_2\leq j$ in the link pattern of $\sigma$.
\end{theorem}
The theorem thus gives a combinatorial characterization of the $B_n$-orbits in $\mathfrak{n}_n^{(2)}$ and their orbit closures in terms of link patterns.


\subsection{Representations of algebras}\label{roa}

As we make key use of results from the representation theory of finite-dimen\-sion\-al algebras for the study of the action of $B_n$ on $\N_n^{(2)}$, we now recall the basic setup of this theory and refer to \cite{ASS} and \cite{AR} for a thorough treatment.
Let $\Q$ be a finite quiver, that is, a directed graph $\Q=(\Q_0,\Q_1,s,t)$ consisting of a finite set of vertices $\Q_0$ and a finite set of arrows $\Q_1$, whose elements are written as $\alpha:s(\alpha)\rar t(\alpha)$; the vertices $s(\alpha)$ and $t(\alpha)$ are called the source and the target of $\alpha$, respectively. A path in $\Q$ is a sequence of arrows $\omega=\alpha_s\ldots\alpha_1$ such that $t(\alpha_{k})=s(\alpha_{k+1})$ for all $k=1,\ldots,s-1$; we formally include a path $\varepsilon_i$ of length zero for each $i\in \Q_0$ starting and ending in $i$. We have an obvious notion of concatenation $\omega\omega'$ of paths $\omega=\alpha_s\ldots\alpha_1$ and $\omega'=\beta_t\ldots\beta_1$ such that $t(\beta_t)=s(\alpha_1)$.\\[1ex]
The path algebra $k\Q$ is defined as the $k$-vector space with basis consisting of all paths in $\Q$, and with multiplication
$$\omega\cdot\omega'=\left\{\begin{array}{cc}\omega\omega'&\mbox{ if $t(\beta_t)=s(\alpha_1)$,}\\
0&\mbox{ otherwise.}\end{array}\right.$$
The radical ${\rm rad}(k\Q)$ is defined as the (two-sided) ideal generated by paths of positive length. An ideal $I$ of $k\Q$ is called admissible if ${\rm rad}(k\Q)^s\subset I\subset{\rm rad}(k\Q)^2$ for some $s$.\\[1ex]
The key feature of such pairs $(\Q,I)$ consisting of a quiver $\Q$ and an admissible ideal $I\subset k\Q$ is the following: every finite-dimensional $k$-algebra $A$ is Morita-equivalent to an algebra of the form $k\Q/I$, in the sense that their categories of finite-dimensional $k$-representations are ($k$-linearly) equivalent.\\[1ex]
A finite-dimensional $k$-representation $M$ of $\Q$ consists of a tuple of $k$-vector spaces $M_i$ for $i\in\Q_0$, and a tuple of $k$-linear maps $M_\alpha:M_i\rar M_j$ indexed by the arrows $\alpha:i\rar j$ in $\Q_1$. A morphism of two such representations $M=((M_i)_{i\in\Q_0},(M_\alpha)_{\alpha\in\Q_1})$ and $N=((N_i)_{i\in\Q_0},(N_\alpha)_{\alpha\in\Q_1})$ consists of a tuple of $k$-linear maps $(f_i:M_i\rar N_i)_{i\in \Q_0}$ such that
$$f_jM_\alpha=N_\alpha f_i\mbox{ for all }\alpha:i\rar j\mbox{ in }\Q_1.$$
For a representation $M$ and a path $\omega$ in $\Q$ as above, we denote $M_\omega=M_{\alpha_s}\cdot\ldots\cdot M_{\alpha_1}$. We call $M$ bound by $I$ if $\sum_\omega\lambda_\omega M_\omega=0$ whenever $\sum_\omega\lambda_\omega\omega\in I$.\\[1ex]
The abelian $k$-linear category of all representations of $\Q$ bound by $I$ is denoted by $\reps$; it is equivalent to the category of finite-dimensional representations of the algebra $k\Q/I$. We have thus found a ``linear algebra model'' for the category of finite-dimensional representations of an arbitrary finite-dimensional $k$-algebra $A$.\\[1ex]
We define the dimension vector $\dimv M\in{\bf N}\Q_0$ of $M$ by $(\dimv M)_{i}=\dim_kM_i$ for $i\in\Q_0$. For a fixed dimension vector $\df\in{\bf N}\Q_0$, we consider the affine space $R_{\df}(Q)= \bigoplus_{\alpha:i\rar j}\Hom_k(k^{d_i},k^{d_j})$; its points naturally correspond to representations $M$ of $\Q$ of dimension vector $\df$ with $M_i=k^{d_i}$ for $i\in \Q_0$. Via this correspondence, the set of such representations bound by $I$ corresponds to a closed subvariety $R_\df(Q,I)\subset R_\df(Q)$.
It is obvious that the algebraic group $\GL_\df=\prod_{i\in\Q_0}\GL_k(k^{d_i})$ acts on $R_{\df}(Q)$ and on $R_\df(Q,I)$ via base change $(g_i)_i\cdot(M_\alpha)_\alpha=(g_jM_\alpha g_i^{-1})_{\alpha:i\rar j}$. By definition, the $\GL_{\df}$-orbits $\Orb_M$ of this action naturally correspond to the isomorphism classes of representations $M$ in $\reps$ of dimension vector $\df$.\\
By the Krull-Schmidt theorem, every representation in $\reps$ is isomorphic to a direct sum of indecomposables, unique up to isomorphisms and permutations. Thus, knowing the isomorphism classes of indecomposable representations in $\reps$ and their dimension vectors, we can classify the orbits of $\GL_\df$ in $R_\df(\Q,I)$.\\[1ex]
For certain classes of finite-dimensional algebras, a convenient tool for the classification of the indecomposable representations is the Auslander-Reiten quiver $\Gamma(\Q,I)$ of $k\Q/I$. Its vertices $[X]$ are given by the isomorphism classes of indecomposable representations of $k\Q/I$; the arrows between two such vertices $[X]$ and $[Y]$ are parametrized by a basis of the space of so-called irreducible maps $f:X\rar Y$. Several standard techniques are available for the calculation of $\Gamma(\Q,I)$, see for example \cite{ASS} and \cite{MRXXXXXX}. We will illustrate one of these techniques, namely the use of covering quivers, in subsection \ref{covsection} in a situation relevant for our setup.

\section{Classification of orbits}\label{Class}

\subsection{Translation to a representation-theoretic problem}\label{trans}

Our aim in this section is to translate the classification problem for the action of $B_n$ on $\N_n^{(2)}$ into a representation-theoretic one. The following is a well-known fact on associated fibre bundles:

\begin{theorem}\label{basicthm}
Let $\G$ be an algebraic group, let $X$ and $Y$ be $\G-$varieties, and let $\pi : X \rar Y$ be a $\G$-equivariant morphism. Assume that $Y$ is a single $\G$-orbit, $Y = \G y_0$. Define $H := \Iso_{\G} (y_0) = \{g \in \G \mid g\cdot y_0 = y_0 \}$ and $F := \pi^{-1} (y_0)$. Then $X$ is isomorphic to the associated fibre bundle $G\times^HF$, and the embedding $\phi: F \hookrightarrow X$ induces a bijection between $H$-orbits in $F$ and $G$-orbits in $X$ preserving orbit closures.
\end{theorem}

We consider the following quiver, denoted by $\Q$ from now on,
\begin{center}
\begin{tikzpicture}
\matrix (m) [matrix of math nodes, row sep=0.05em,
column sep=2em, text height=1.5ex, text depth=0.2ex]
{\Q: & \bullet & \bullet &  \bullet & \cdots  & \bullet & \bullet  & \bullet \\ & { 1} & { 2} &  { 3} & &   { n-2} &  { n-1}  & { n} \\ };
\path[->]
(m-1-2) edge node[above=0.05cm] {$\alpha_1$} (m-1-3)
(m-1-3) edge  node[above=0.05cm] {$\alpha_2$}(m-1-4)
(m-1-6) edge  node[above=0.05cm] {$\alpha_{n-2}$}(m-1-7)
(m-1-7) edge node[above=0.05cm] {$\alpha_{n-1}$} (m-1-8)
(m-1-8) edge [loop right] node{$\alpha$} (m-1-8);
\end{tikzpicture}
\end{center}
together with the ideal $I\subset k\Q$ generated by the path $\alpha^2$. We consider the full subcategory $\repinj$ of $\reps$ consisting of representations $M$ for which the linear maps $M_{\alpha_1},\ldots,M_{\alpha_{n-1}}$ are injective. Corresponding to this subcategory, we have an open subset $R_\df^{\rm inj}(\Q,I)\subset R_\df(\Q,I)$, which is stable under the $\GL_\df$-action. We consider the dimension vector $\df:=\vect\in {\bf N}^n$.
\begin{lemma} \label{thm}
There exists a closure-preserving bijection between the set of $B_n$-orbits in $\N_n^{(2)}$ and the set of $\GL_\df$-orbits in $R_\df^{\rm inj}(Q,I)$.
\end{lemma}
\begin{proof}
Consider the subquiver $\Qtil$ of $\Q$ with $\Qtil_0=\Q_0$ and $\Qtil_1=\Q_1\setminus\{\alpha\}$. We have a natural $\GL_\df$-equivariant projection $\pi:R_\df^{\rm inj}(Q,I)\rightarrow R_{{\df}}^{\rm inj}(\Qtil)$. The variety $R_{{\df}}^{\rm inj}(\Qtil)$ consists of tuples of injective maps, thus the action of $\GL_\df$ on $R_{{\df}}(\Qtil)$ is easily seen to be transitive. Namely, $R_{{\df}}^{\rm inj}(\Qtil)$ is the orbit of the representation
$$y_0:=k \xrar{\alpha_{1}} k^{2} \xrar{\alpha_{2}} \cdots \xrar{\alpha_{n-2}} k^{n-1} \xrar{\alpha_{n-1}} k^{n},$$ with $\alpha_i$ being the canonical embedding from $k^i$ to $k^{i+1}$. The stabilizer $H$ of $y_0$ is isomorphic to $B_n$, and the fibre of $\pi$ over $y_0$ is isomorphic to $\N_n^{(2)}$. Thus, $R_\df^{\rm inj}(\Q,I)$ is isomorphic to the associated fibre bundle $\GL_\df\times^{B_n}\N_n^{(2)}$, yielding the claimed bijection.
\end{proof}

\subsection{Classification of indecomposables in $\reps$}\label{covsection}

By the results of the previous section, it suffices to classify the indecomposable representations in $\reps$ to obtain a classification of the orbits of $B_n$ in $\N_n^{(2)}$. We compute the Auslander-Reiten quiver $\Gamma$ of $k\Q/I$ using covering theory, which is described in \cite{MRXXXXXX} as mentioned before. We consider the (infinite) quiver $\widehat{\Q}$ given by
\begin{center}
\begin{tikzpicture}
\matrix (m) [matrix of math nodes, row sep=0.85em,
column sep=2em, text height=1.5ex, text depth=0.2ex]
{&&&&\vdots&&&~ \\  & \bullet & \bullet &  \bullet & \cdots  & \bullet & \bullet  & \bullet \\
 \widehat{\Q}: & \bullet & \bullet &  \bullet & \cdots  & \bullet & \bullet  & \bullet \\
  & \bullet & \bullet &  \bullet & \cdots  & \bullet & \bullet  & \bullet \\
  &&&&\vdots&&&~ \\
 & {\rm 1} & {\rm 2} &  {\rm 3} & &   {\rm n-2} &  {\rm n-1}  & {\rm n} \\ };
\path[->]
(m-2-2) edge node[above=0.05cm] { }(m-2-3)
(m-2-3) edge  node[above=0.05cm] { } (m-2-4)
(m-2-6) edge  node[above=0.05cm] { } (m-2-7)
(m-2-7) edge node[above=0.05cm]  { } (m-2-8)

(m-3-2) edge node[above=0.05cm]  { } (m-3-3)
(m-3-3) edge  node[above=0.05cm] { } (m-3-4)
(m-3-6) edge  node[above=0.05cm] { } (m-3-7)
(m-3-7) edge node[above=0.05cm] { }  (m-3-8)

(m-4-2) edge node[above=0.05cm]  { } (m-4-3)
(m-4-3) edge  node[above=0.05cm] { } (m-4-4)
(m-4-6) edge  node[above=0.05cm] { } (m-4-7)
(m-4-7) edge node[above=0.05cm] { }  (m-4-8)

(m-1-8) edge node[right=0.05cm] {} (m-2-8)
(m-2-8) edge node[right=0.05cm] {$\alpha_i$} (m-3-8)
(m-3-8) edge node[right=0.05cm] {$\alpha_{i+1}$} (m-4-8)
(m-4-8) edge node[right=0.05cm] {} (m-5-8);
\end{tikzpicture}
\end{center}

with the ideal $\widehat{I}$ generated by all paths $\alpha_{i+1}\alpha_i$, and the quiver $\Q'$ given by
\begin{center}
\begin{tikzpicture}
\matrix (m) [matrix of math nodes, row sep=0.85em,
column sep=2em, text height=1.5ex, text depth=0.2ex]
{\Q': & \bullet & \bullet &  \bullet & \cdots  & \bullet & \bullet  & \bullet \\
 & \bullet & \bullet &  \bullet & \cdots  & \bullet & \bullet  & \bullet \\
 & {\rm 1} & {\rm 2} &  {\rm 3} & &   {\rm n-2} &  {\rm n-1}  & {\rm n} \\ };
\path[->]
(m-1-2) edge node[above=0.05cm] { }(m-1-3)
(m-1-3) edge  node[above=0.05cm] { } (m-1-4)
(m-1-6) edge  node[above=0.05cm] { } (m-1-7)
(m-1-7) edge node[above=0.05cm]  { } (m-1-8)

(m-2-2) edge node[above=0.05cm]  { } (m-2-3)
(m-2-3) edge  node[above=0.05cm] { } (m-2-4)
(m-2-6) edge  node[above=0.05cm] { } (m-2-7)
(m-2-7) edge node[above=0.05cm] { }  (m-2-8)
(m-1-8) edge node[right=0.05cm] {$\alpha$} (m-2-8);
\end{tikzpicture}
\end{center}

\begin{figure}[ht]
\includegraphics[width=1\textwidth]{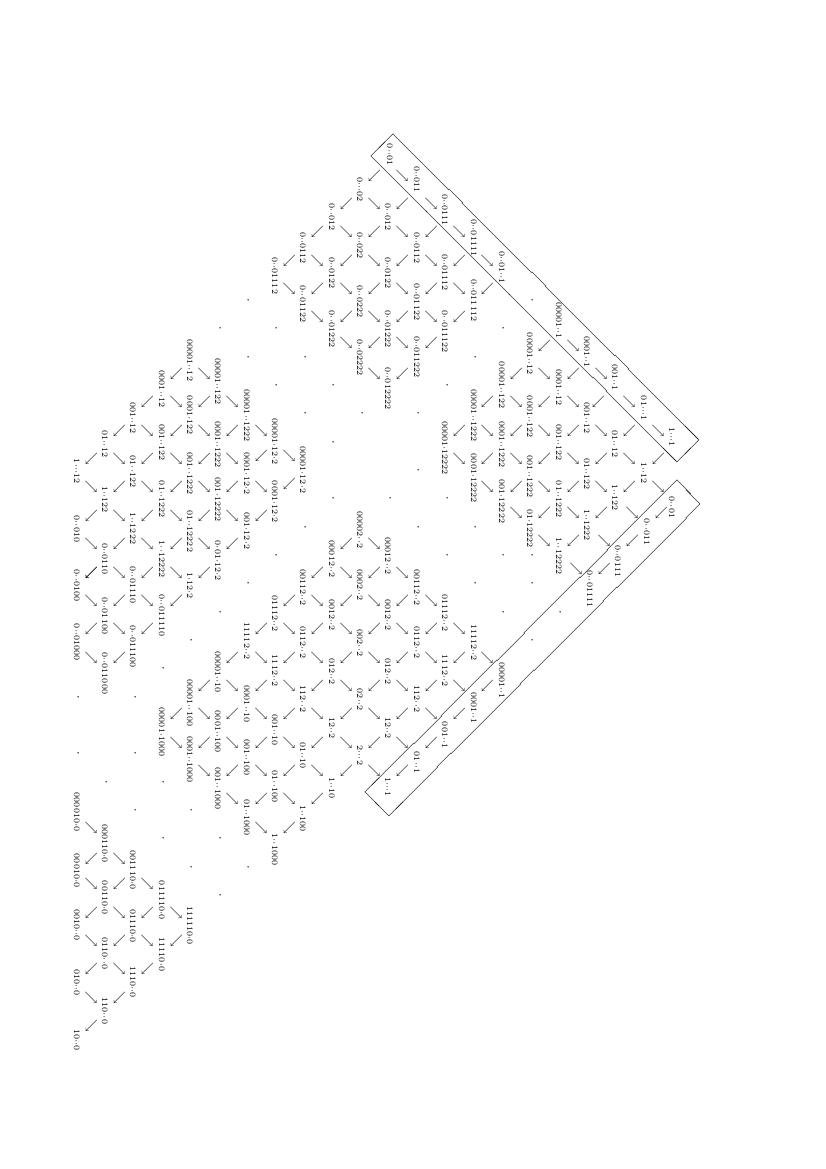}
\caption{\label{D.dttn.zeichnung} The Auslander-Reiten quiver of $\reps$}  \label{ARK}
\end{figure}

The quiver $\widehat{\Q}$ carries a natural action of the group ${\bf Z}$ by shifting the rows, such that $\widehat{\Q}/{\bf Z}\cong\Q$. Moreover, $\Q'$ naturally embeds into $\widehat{\Q}$, such that the composition of this inclusion with the projection $\widehat{\Q}\rightarrow\Q$ is surjective. By results of covering theory \cite{MRXXXXXX}, we have corresponding maps of the Auslander-Reiten quivers, namely an embedding $\Gamma(\Q')\rightarrow\Gamma(\widehat{\Q},\widehat{I})$ and a quotient $\Gamma(\widehat{\Q},\widehat{I})\rightarrow\Gamma(\Q,I)$, such that the composition is surjective.
Since $\Q'$ is nothing else than a Dynkin quiver of type $A_{2n}$, it is routine to calculate its Auslander-Reiten quiver (see \cite{ASS}), and we derive the Auslander-Reiten quiver $\Gamma=\Gamma(\Q,I)$ just by making the identifications resulting from the action of ${\bf Z}$, which can be read off from the dimension vectors of indecomposable representations. More examples and details concerning the calculation of Auslander-Reiten quivers using covering theory can also be found in \cite{MRXXXXXX}.\\
We finally arrive at the picture (the marked regions have to be identified) given in figure \ref{ARK}.



We define the following representations $\U_{i,j}$ for $1\leq i,j\leq n$, $\V_i$ for $1\leq i\leq n$ and $\W_{i,j}$ for $1\leq i\leq j\leq n$ in $\reps$ (graphically represented by dots for basis elements and arrows for a map sending one basis element to another one):\\[1ex]
$\U_{i,j}$ for $1\leq j\leq i\leq n$:
$$\begin{tikzpicture}
\matrix (m) [matrix of math nodes, row sep=0.02em,
column sep=0.08em, text height=1.0ex, text depth=0.25ex]
{ 0 & \xrar{0} & \cdots & \xrar{0} & 0 & \xrar{0} & k & \xrar{id} & \cdots & \xrar{id} & k & \xrar{e_1} & k^{2} & \xrar{id}& \cdots & \xrar{id} & k^{2}  \\
 & & & & & & \bullet & \rar & \cdots & \rar & \bullet & \rar & \bullet&  \rar & \cdots & \rar& \bullet \\
  & & & & & & j & & & & &  & i & & & & n\\
 & & & & & & & & & && & \bullet & \rar & \cdots & \rar &  \bullet \\ };
\path[->]
(m-2-17) edge [bend left=80] (m-4-17)
(m-1-17) edge [loop right] node{$\alpha$} (m-1-17);
\end{tikzpicture}$$
$\U_{i,j}$ for $1\leq i< j\leq n$:
$$\begin{tikzpicture}
\matrix (m) [matrix of math nodes, row sep=0.02em,
column sep=0.08em, text height=1.0ex, text depth=0.25ex]
{  0 & \xrar{0} & \cdots & \xrar{0} & 0 & \xrar{0} & k & \xrar{id} & \cdots & \xrar{id} & k & \xrar{e_2} & k^{2} & \xrar{id}& \cdots & \xrar{id} & k^{2}  \\
 & & & & & & \bullet & \rar & \cdots & \rar & \bullet & \rar & \bullet&  \rar & \cdots & \rar& \bullet \\
  & & & & & & i & & & & &  & j & & & & n\\
 & & & & & & & & & && & \bullet & \rar & \cdots & \rar &  \bullet \\ };
\path[->]
(m-4-17) edge [bend right=80] (m-2-17)
(m-1-17) edge [loop right] node{$\alpha$} (m-1-17);
\end{tikzpicture}$$
$\V_{i}$ for $1\leq i\leq n$:
$$\begin{tikzpicture}
\matrix (m) [matrix of math nodes, row sep=0.02em,
column sep=0.1em, text height=1.0ex, text depth=0.25ex]
{  0 & \xrar{0} & \cdots & \xrar{0} & 0 & \xrar{0} & k & \xrar{id} & \cdots & \xrar{id} & k  \\
  & & & & & & i & & & &  n\\
 & & & & & & \bullet & \rar & \cdots & \rar & \bullet  \\ };
\path[->]
(m-1-11) edge [loop right] node{$0$} (m-1-11);
\end{tikzpicture}$$
$\W_{i,j}$ for $1\leq i\leq j<n$:
$$\begin{tikzpicture}
\matrix (m) [matrix of math nodes, row sep=0.02em,
column sep=0.1em, text height=1.0ex, text depth=0.25ex]
{ 0 & \xrar{0} & \cdots & \xrar{0} & 0 & \xrar{0} & k & \xrar{id} & \cdots & \xrar{id} & k & \xrar{0} & 0 & \xrar{0} & \cdots & \xrar{0} & 0  \\
  & & & & & & i & & & & j & & & & & &  n\\
 & & & & & & \bullet & \rar & \cdots & \rar & \bullet &&&&&&  \\ };
\path[->]
(m-1-17) edge [loop right] node{$0$} (m-1-17);
\end{tikzpicture}$$
Here we denote $e_1=
\left(\begin{array}{l}
1 \\
0
                                                                                                                                                         \end{array}\right)$, $e_2=
\left(\begin{array}{l}
0 \\
1
                                                                                                                                                         \end{array}\right)$ and $\alpha=
\left(
\begin{array}{ll}
0 & 0 \\
1 & 0 \end{array}\right)$.

\begin{theorem}\label{UVW} The representations $\U_{i,j}$, $\V_i$ and $\W_{i,j}$ form a system of representatives of the indecomposable objects in $\reps$. The representations $\U_{i,j}$ and $\V_i$ form a system of representatives of the indecomposable objects in $\repinj$
\end{theorem}

\begin{proof} The endomorphism rings of these representations are easily computed to be
$${\rm End}(\U_{i,j})\cong k\mbox{ for }i>j,\;\;\; {\rm End}(\U_{i,j})\cong k[x]/(x^2)\mbox{ for }i\leq j,$$
$${\rm End}(\V_i)\cong k,\;\;\; {\rm End}(\W_{i,j})\cong k,$$
thus they are indecomposable. Their dimension vectors are
$$(0\ldots 01\ldots 12\ldots 2),\;(0\ldots\ldots 01\ldots\ldots 1)\mbox{ and }(0\ldots\ldots 01\ldots\ldots 10\ldots\ldots 0),$$ respectively. These are precisely the dimension vectors appearing in $\Gamma(\Q,I)$, thus we have found all indecomposables. It is clear from the definition that the indecomposable representations belonging to $\repinj$ are the $\U_{i,j}$ and the $\V_i$.\end{proof}


\subsection{Classification of $B_n$-orbits in $\N_n^{(2)}$}\label{cob}

Our next aim is to parametrize the isomorphism classes of representations in $\repinj$ of dimension vector $\df$. As mentioned before, the Krull-Schmidt theorem states that every representation can be decomposed into a direct sum of indecomposables in an essentially unique way.

\begin{theorem}\label{classification} The isomorphism classes $M$ in $\repinj$ of dimension vector $\df$ are in natural bijection to
\begin{enumerate}
\item $n\times n$-matrices $A=(m_{i,j})_{i,j}$ with entries $0$ or $1$, such that $\sum_jm_{i,j}+\sum_jm_{j,i}\leq 1$ for all $i=1,\ldots,n$,
\item oriented link patterns on $\{1,\ldots,n\}$, that is, oriented graphs on the set $\{1,\ldots,n\}$ such that every vertex is incident with at most one arrow.
\end{enumerate}
Moreover, if an isomorphism class $M$ corresponds to a matrix $A$ under this bijection, the orbit $\Orb_M\subset R_{\df}^{\rm inj}(\Q,I)$ and the orbit $\Orb_A\in\N_n^{(2)}$ correspond to each other via the bijection of Lemma \ref{thm}.
\end{theorem}

\begin{proof} Let $M$ be a representation in $\repinj$ of dimension vector $\df$, so $$M=\bigoplus_{i,j=1}^n\U_{i,j}^{m_{i,j}}\oplus\bigoplus_{i=1}^n\V_i^{n_i}$$ for some multiplicities $m_{i,j},n_i\in\bf N$ by Theorem \ref{UVW}. Since $\dimv M=\vect$, we simply need to calculate all tuples $(m_{i,j},n_{i})$ such that
$$\sum_{i,j=1}^{n} m_{i,j}\dimv\U_{i,j}~+~\sum_{i=1}^{n}
n_{i}\dimv\V_{i}~=~\df~=~\vect.$$
Applying the automorphism $\delta$ of ${\bf Z}^n$ defined by
$$ \delta(d_{1},d_{2},\ldots,d_{n})=
(d_{1},d_{2}-d_{1},d_{3}-d_{2},\ldots,d_{n}-d_{n-1}),$$ this condition is
equivalent to
$$\sum_{i,j=1}^{n} m_{i,j}\delta(\dimv\U_{i,j}) +
\sum_{i=1}^{n} n_{i}\delta(\dimv\V_{i})=(1,1,\ldots,1,1).$$
If we fix $i \in \{1,\ldots,n\}$, this condition states that
$$1~=~\sum_{j=1}^{n} m_{i,j}~+~\sum_{j=1}^{n} m_{j,i}~+~n_{i}.$$
We can extract an oriented graph on the set of vertices $\{1,\ldots,n\}$
from $(m_{i,j})_{i,j}$ as follows: for all $1\leq i,j\leq n$, we have an
arrow from $j$ to $i$ if $m_{i,j}=1$. The conditions on $(m_{i,j})_{i,j}$
ensure that this graph is in fact an oriented link pattern. The matrix
$(m_{i,j})_{i,j}$ is obviously $2$-nilpotent.\\[1ex]
The decomposition of $M$ into indecomposables can be visualized as follows.
$$\begin{tikzpicture}
\matrix (m) [matrix of math nodes, row sep=0.25em,
column sep=0.5em, text height=1.5ex, text depth=0.25ex]
{& \bullet & \rightarrow & \bullet& \rightarrow & \bullet & \cdots & \bullet & \rightarrow & \bullet & \rightarrow & \bullet & ~~ & 1  \\
& &  & \bullet& \rightarrow & \bullet & \cdots & \bullet & \rightarrow & \bullet & \rightarrow & \bullet & ~~ & 2\\
& &  & & & \bullet & \cdots & \bullet & \rightarrow & \bullet & \rightarrow & \bullet & ~~ & 3\\
M:& & & &  &  &  &  & \vdots  &  & \vdots  & & & \\
& &  & &  &  &  & \bullet & \rightarrow & \bullet & \rightarrow & \bullet & ~~ & n-2\\
& &  & &  &  &  &  &  & \bullet & \rightarrow & \bullet & ~~ & n-1\\
& &  & &  & &  &  &  &  &  & \bullet & ~~ & n\\ };
\path[->]
(m-1-12) edge [bend left=80] (m-3-12)
(m-7-12) edge [bend right=50] (m-2-12)
(m-6-12) edge [bend right=120] (m-5-12);
\end{tikzpicture}$$
The arrows in the rightmost column of the diagram allow us to read off the indecomposable direct summands of $M$. Namely, $\U_{i,j}$ is a direct summand of $M$ if and only if there is an arrow $j\rar i$. If there is no arrow at $k$, the indecomposable $\V_k$ is a direct summand of $M$.\\[1ex]
Shortening the above picture to the rightmost column, $M$ corresponds to an oriented link pattern:
$$\xygraph{ !{<0cm,0cm>;<1cm,0cm>:<0cm,1cm>::}
!{(-3,0) }*+{\underset{1}\bullet}="1"
!{(-2,0)}*+{\underset{2}\bullet}="2"
!{(-1,0)}*+{\underset{3}\bullet}="3"
!{(0,0) }*+{\ldots}="4"
!{(1,0) }*+{\underset{n-2}\bullet}="5"
!{(2,0)}*+{\underset{n-1}\bullet}="6"
!{(3,0)}*+{\underset{n}\bullet}="7"
"1":@/^0.75cm/"3"
"7":@/_1cm/"2"
"6":@/_0.65cm/"5"
}$$
\end{proof}


For a given matrix $A\in\N_n^{(2)}$, we would like to decide to which oriented link pattern it corresponds.
Define $U_i=\la e_1,\ldots,e_i\ra$, the span of the first $i$ coordinate vectors in $k^n$, and define a matrix $D^A =(d_{i,j}^A)_{i,j}$ by setting $d_{i,j}^A:=\dim(U_i\cap A(U_j))$ (we formally define $d_{i,j}^A=0$ for $i=0$ or $j=0$). The matrix $D^A$ is obviously an invariant for the $B_n$-action on $\N_n^{(2)}$. It is easy to extract an oriented link pattern from $D^A$ as follows:

\begin{lemma} The matrix $A$ belongs to the orbit of a matrix $(m_{i,j})_{i,j}$ as above if and only if $d_{i,j}^A=\sum_{i'\leq i;\, j'\leq j}m_{i',j'}$ or, conversely, $m_{i,j}=d_{i,j}^A-d_{i-1,j}^A-d_{i,j-1}^A+d_{i-1,j-1}^A$ for all $1\leq i,j\leq n$.
\end{lemma}

\begin{proof} By $B_n$-invariance, we just have to compute $D^A$ for $A=(m_{i,j})_{i,j}$ as in the previous theorem. We have $e_k\in U_i\cap A(U_j)$ if and only $k\leq i$ and there exists $l\leq j$ such that $m_{k,l}=1$ or, equivalently, such that there exists an arrow $l\rar k$ in the corresponding oriented link pattern. Since both $U_i$ and $A(U_j)$ are spanned by coordinate vectors $e_k$, we thus have $d_{i,j}^A=\dim(U_i\cap A(U_j))=\sum_{k\leq i,\, l\leq j}m_{k,l}$. The second formula follows.
\end{proof}

We can also rederive Theorem \ref{mel1} of A. Melnikov: every $B_n-$orbit of an upper-triangular $2$-nilpotent matrix corresponds to the orbit of a representation in $\repinj$ of dimension vector $\df$ which does not contain $\U_{i,j}$ for $i\geq j$ as a direct summand. In this case, the corresponding link pattern consists of arrows pointing in the same direction. We can thus delete the orientation and arrive at a link pattern as in \cite{Melnikov1}.\\[2ex]
\remark Our method is easily generalized to obtain a classification of orbits for a
more general group action: let $\Pa\subset \GL_n$ be the parabolic subgroup consisting of block-upper
triangular matrices with block-sizes $(b_1,\ldots,b_k)$.
Then $\Pa$ acts on $\N_n^{(2)}$ by conjugation, and the same reasoning as
above yields a bijection between $\Pa$-orbits in $\N_n^{(2)}$ and
isomorphism classes of representations in $\repinj$ of dimension vector
$(b_1,b_1+b_2,\ldots,\sum_{i=1}^k b_i)$.
Using the analysis of this section, the $\Pa$-orbits in $\N_n^{(2)}$
correspond bijectively to matrices $(m_{i,j})_{i,j}$ such that
$$\sum_j m_{i,j}+\sum_j m_{j,i} \leq b_i$$
for all $i=1,\ldots,k$. Consequently, they correspond bijectively to ``enhanced oriented link
patterns of type $(b_1,\ldots,b_k)$'', namely, to oriented graphs on the set
$\{1,\ldots,k\}$ such that the vertex $i$ is incident with at most $b_i$
arrows for all $i$.


\section{Orbit closures}\label{close}

After classifying the orbits via oriented link patterns, we describe the corresponding orbit closures. Again, we will solve this problem using results about the geometry of representations of algebras. Two theorems of G. Zwara are the key to calculating these orbit closures, see \cite{Zwara1} and \cite{Zwara2} for more details.

\subsection{A criterion for degenerations}\label{cfd}

Let $M$ and $M'$ be two representations in $\reps$ of the same dimension vector $\df$. We say that $M$ degenerates to $M'$ if $\Orb_{M'}\subset\overline{\Orb_M}$ in $R_{\df}(\Q,I)$, which will be denoted by $M\leq_{\rm deg}M'$. Since the correspondence of Lemma \ref{thm} preserves orbit closure relations, we know that $M\leq_{\rm deg}M'$ if and only if the corresponding $2$-nilpotent matrices, denoted by $A=(m_{i,j})_{i,j}$ and $A'=(m'_{i,j})_{i,j}$, respectively, fulfill $\Orb_{A'}\subset\overline{\Orb_{A}}$ in $\N_n^{(2)}$.

\begin{theorem}\label{zw}(Zwara)
Suppose an algebra $k\Q/I$ is representation-finite, that is, $k\Q/I$ admits only finitely many isomorphism classes of indecomposable representations. Let $M$ and $M'$ be two finite-dimensional representations of $k\Q/I$ of the same dimension vector.\\
Then $M\leq_{\rm deg} M'$ if and only if $\dim_k \Hom(U, M ) \leq \dim_k \Hom(U, M')$ for every representation $U$ of $k\Q/I$.
\end{theorem}

To simplify notation, we set $[U, V] := \dim_k \Hom(U, V )$ for two representations $U$ and $V$. Since the dimension of a homomorphism space is additive with respect to direct sums, we only have to consider the inequality $[U,M]\leq[U,M']$ for indecomposable representations $U$ to characterize a degeneration $M\leq_{\rm deg}M'$. Furthermore, since $[\W_{i,j},M]=0$ for all representations $M$ in $\repinj$ by a direct calculation, we can restrict these indecomposables $U$ to those of type $\U_{i,j}$ and $\V_i$ of the previous section.\\[1ex]
We can easily calculate the dimensions of homomorphism spaces between these indecomposable representations.
\begin{lemma}\label{hom} For $i, j, k, l \in \{1,\ldots, n\}$ we have
\begin{itemize}
\item $[\V_k , \V_i ] = \delta_{i\leq k}$,
\item $[\V_k , \U_{i,j} ] = \delta_{i\leq k}$,
\item $[\U_{k,l} , \V_i ] = \delta_{i\leq l}$,
\item $[\U_{k,l} , \U_{i,j} ] = \delta_{i\leq l} + \delta_{j\leq l} \cdot \delta_{i\leq k}$,
\end{itemize}
where
$\delta_{x\leq y} := \left\{ \begin{array}{ll} 1, & \hbox{${\rm if}~ x\leq y$;} \\ 0, & \hbox{${\rm otherwise}$.} \end{array} \right. $
\end{lemma}

For a representation $M$ in $\repinj$ of dimension vector $\df$ (or equivalently, for the corresponding $2$-nilpotent matrix $A$), consider the corresponding oriented link pattern. Define $p_k^M$ as the number of vertices to the left of $k$ which are not incident with an arrow, plus the number of arrows whose target vertex is to the left of $k$. Define $q^M_{k,l}$ as $p_l^M$ plus the number of arrows whose source vertex lies to the left of $l$ and whose target vertex lies to the left of $k$.

\begin{theorem}\label{pq} We have $M\leq_{\rm deg} M'$ (or equivalently, $\Orb_{A'}\subset\overline{\Orb_{A}}$ in the notation above) if and only if $p_k^M\leq p_k^{M'}$ and $q_{k,l}^M\leq q_{k,l}^{M'}$ for all $k,l=1,\ldots,n$.
\end{theorem}

\begin{proof} Given two representations $M$ and $M'$, we write
$$M=\bigoplus_{i,j=1}^n\U_{i,j}^{m_{i,j}}\oplus\bigoplus_{i=1}^n\V_i^{n_i}\mbox{ and }M'=\bigoplus_{i,j=1}^n\U_{i,j}^{m'_{i,j}}\oplus\bigoplus_{i=1}^n\V_i^{n'_i}.$$
For an indecomposable $U$, the condition $[U, M ] \leq [U, M' ]$ is then equivalent to
$$\sum_{i,j=1}^{n} m_{i,j} [U, \U_{i,j}]+ \sum_{i=1}^n n_i[U, \V_i ] \leq \sum_{i,j=1}^{n} m'_{i,j} [U, \U_{i,j}]+ \sum_{i=1}^n n'_i[U, \V_i ].$$
Using the dimensions of homomorphism spaces between indecomposable representations stated in the previous lemma, we calculate
$$p_k^M=\sum_{i\leq k; j}m_{i,j}+\sum_{i\leq k}n_i$$
and
$$q_{k,l}(M)=p_l^M+\sum_{i\leq k; j\leq l}m_{i,j},$$
and the condition $[U,M]\leq [U,M']$ is equivalent to the conditions $p_k^M\leq p_k^{M'}$ and $q_{k,l}^M\leq q_{k,l}^{M'}$ for all $k,l=1,\ldots,n$. The claimed interpretation of these values $p_k^M$ and $q_{k,l}^M$ in terms of oriented link patterns follows from
Theorem \ref{classification}.
\end{proof}


\subsection{Minimal degenerations}\label{md}

As a next step, we develop a combinatorial method to produce all degenerations of a given representation $M$ in $\repinj$ of dimension vector $\df$ out of its oriented link pattern. It is sufficient to construct all minimal degenerations, that is, degenerations $M<_{\rm deg}M'$ such that if $M\leq_{\rm deg}L\leq_{\rm deg}M'$, then $M\cong L$ or $M'\cong L$. Minimal degenerations are denoted by $M<_{\rm mdeg}M'$.\\[1ex]
In \cite{Zwara2}, G. Zwara describes all minimal degenerations; the result is stated here in a generality sufficient for our purposes.
Denote by $\leq_{\rm ext}$ the transitive closure of the relation on representations given by $M\leq M'$ if there exists a short exact sequence $0\rar M_1'\rar M\rar M_2'\rar 0$ such that $M'\cong M_1'\oplus M_2'$.
\begin{theorem}
Let $M$ and $M'$ be representations in $\repinj$.\\
If $M<_{\rm mdeg}M' $, then one of the following holds:
\begin{enumerate}
 \item $M<_{ext} M'$
 \item There are representations $W$, $\widetilde{M}$, $\widetilde{M}'$ in $\repinj$ such that
\begin{enumerate}
\item $M\cong W\oplus \widetilde{M}$
\item $M'\cong W\oplus \widetilde{M'}$
\item $\widetilde{M}<_{\rm mdeg} \widetilde{M'}$
\item $\widetilde{M'}$ is indecomposable.
\end{enumerate}
\end{enumerate}
\end{theorem}

Combining this theorem with the technique of \cite[Theorem 4]{Bongartz}, we obtain a characterization of minimal disjoint degenerations, that is, minimal degenerations $M<_{\rm mdeg}M'$ such that $M$ and $M'$ do not share a common direct summand:

\begin{corollary}
Let $M<_{\rm mdeg}M'$ be a minimal disjoint degeneration as before. Then either $M'$ is indecomposable or $M'\cong U\oplus V$, where $U$ and $V$ are indecomposables and there exists an exact sequence $0\rar U\rar M\rar V\rar 0$ or $0\rar V\rar M\rar U\rar 0$.
\end{corollary}

Thus we see that all minimal degenerations are of the form $W\oplus M<_{\rm mdeg}W\oplus M'$, where $M$ and $M'$ are as in the corollary, thus $M'$ involves at most two indecomposable direct summands. Translating this to the language of oriented link patterns using Theorem \ref{classification}, we have ``localized'' the problem to the consideration of at most four vertices of an oriented link pattern. In this local case, we can apply Theorem \ref{pq} and easily work out all minimal degenerations.

\begin{theorem}
Every minimal degeneration is of the form given in one of the following diagrams showing parts of the degeneration posets in terms of oriented link patterns. We assume that $a<b$ (resp. $a<b<c$, resp. $a<b<c<d$) are vertices of an oriented link pattern, and only indicate the changes to the arrows incident with one of these vertices; all other arrows are left unchanged. 

$$
\begin{xy}
(25,0)*{\xygraph{ !{<0cm,0cm>;<1cm,0cm>:<0cm,1cm>::}
!{(-0.5,0)}*+{\underset{a}\bullet}="1"
!{(0.5,0)}*+{\underset{b}\bullet}="2"
"1":@/^0.6cm/"2"}};
(25,-30)*{\xygraph{ !{<0cm,0cm>;<1cm,0cm>:<0cm,1cm>::}
!{(-0.5,0) }*+{\underset{a}\bullet}="1"
!{(0.5,0)}*+{\underset{b}\bullet}="2"
"2":@/_0.6cm/"1"
 }};
(25,-60)*{\xygraph{ !{<0cm,0cm>;<1cm,0cm>:<0cm,1cm>::}
!{(-0.5,0) }*+{\underset{a}\bullet}="1"
!{(0.5,0)}*+{\underset{b}\bullet}="2"
 }};
(85,0)*{\xygraph{ !{<0cm,0cm>;<1cm,0cm>:<0cm,1cm>::}
!{(-0.6,0) }*+{\underset{a}\bullet}="1"
!{(0,0) }*+{\underset{b}\bullet}="2"
!{(0.6,0) }*+{\underset{c}\bullet}="3"
"1":@/^0.7cm/"3"
 }};
(65,-20)*{\xygraph{ !{<0cm,0cm>;<1cm,0cm>:<0cm,1cm>::}
!{(-0.6,0)}*+{\underset{a}\bullet}="1"
!{(0,0)}*+{\underset{b}\bullet}="2"
!{(0.6,0)}*+{\underset{c}\bullet}="3"
"1":@/^0.7cm/"2"
 }};
(+105,-20)*{\xygraph{ !{<0cm,0cm>;<1cm,0cm>:<0cm,1cm>::}
!{(-0.6,0) }*+{\underset{a}\bullet}="1"
!{(0,0) }*+{\underset{b}\bullet}="2"
!{(0.6,0) }*+{\underset{c}\bullet}="3"
"2":@/^0.7cm/"3"
 } };
(65,-40)*{\xygraph{ !{<0cm,0cm>;<1cm,0cm>:<0cm,1cm>::}
!{(-0.6,0) }*+{\underset{a}\bullet}="1"
!{(0,0) }*+{\underset{b}\bullet}="2"
!{(0.6,0) }*+{\underset{c}\bullet}="3"
"2":@/_0.7cm/"1"
 }};
(+105,-40)*{\xygraph{ !{<0cm,0cm>;<1cm,0cm>:<0cm,1cm>::}
!{(-0.6,0) }*+{\underset{a}\bullet}="1"
!{(0,0) }*+{\underset{b}\bullet}="2"
!{(0.6,0) }*+{\underset{c}\bullet}="3"
"3":@/_0.7cm/"2"
 }};
(85,-60)*{\xygraph{ !{<0cm,0cm>;<1cm,0cm>:<0cm,1cm>::}
!{(-0.6,0) }*+{\underset{a}\bullet}="1"
!{(0,0) }*+{\underset{b}\bullet}="2"
!{(0.6,0) }*+{\underset{c}\bullet}="3"
"3":@/_0.7cm/"1"
 }};
{\ar@{-}@/_0pc/(25,-5);(25,-22)}
{\ar@{-}@/_0pc/(25,-35);(25,-52)}
{\ar@{-}@/_0pc/(79,-3);(66,-13.5)}
{\ar@{-}@/_0pc/(86.5,-3);(100,-13.5)}
{\ar@{-}@/_0pc/(61.0,-23.5);(61,-32.5)}
{\ar@{-}@/_0pc/(+100.5,-23.5);(100.5,-32.5)}
{\ar@{-}@/_0pc/(69,-22.5);(+90,-34.5)}
{\ar@{-}@/_0pc/(90,-22.5);(69,-34.5)}
{\ar@{-}@/_0pc/(60,-43.5);(72,-52.5)}
{\ar@{-}@/_0pc/(94,-43.5);(80.5,-52.5)}
\end{xy}
 $$

$$\begin{xy}
(0,0)*{\xygraph{ !{<0cm,0cm>;<1cm,0cm>:<0cm,1cm>::}
!{(-0.8,0) }*+{\underset{a}\bullet}="1"
!{(-0.275,0)}*+{\underset{b}\bullet}="2"
!{(0.275,0)}*+{\underset{c} \bullet}="3"
!{(0.8,0) }*+{\underset{d}\bullet}="4"
"1":@/^0.7cm/"4"
"2":@/^0.5cm/"3"
 }};
(-20,-20)*{\xygraph{ !{<0cm,0cm>;<1cm,0cm>:<0cm,1cm>::}
!{(-0.8,0) }*+{\underset{a}\bullet}="1"
!{(-0.275,0)}*+{\underset{b}\bullet}="2"
!{(0.275,0)}*+{\underset{c} \bullet}="3"
!{(0.8,0) }*+{\underset{d}\bullet}="4"
"1":@/^0.7cm/"3"
"2":@/^0.7cm/"4"
 }};
(+20,-20)*{\xygraph{ !{<0cm,0cm>;<1cm,0cm>:<0cm,1cm>::}
!{(-0.8,0) }*+{\underset{a}\bullet}="1"
!{(-0.275,0)}*+{\underset{b}\bullet}="2"
!{(0.275,0)}*+{\underset{c} \bullet}="3"
!{(0.8,0) }*+{\underset{d}\bullet}="4"
"1":@/^0.7cm/"4"
"3":@/_0.5cm/"2"
 }};
(-30,-40)*{\xygraph{ !{<0cm,0cm>;<1cm,0cm>:<0cm,1cm>::}
!{(-0.8,0) }*+{\underset{a}\bullet}="1"
!{(-0.275,0)}*+{\underset{b}\bullet}="2"
!{(0.275,0)}*+{\underset{c} \bullet}="3"
!{(0.8,0) }*+{\underset{d}\bullet}="4"
"2":@/^0.7cm/"4"
"3":@/_0.7cm/"1"
 }};
(+0,-40)*{\xygraph{ !{<0cm,0cm>;<1cm,0cm>:<0cm,1cm>::}
!{(-0.8,0) }*+{\underset{a}\bullet}="1"
!{(-0.275,0)}*+{\underset{b}\bullet}="2"
!{(0.275,0)}*+{\underset{c} \bullet}="3"
!{(0.8,0) }*+{\underset{d}\bullet}="4"
"1":@/^0.7cm/"2"
"3":@/^0.7cm/"4"
 }};
(+30,-40)*{\xygraph{ !{<0cm,0cm>;<1cm,0cm>:<0cm,1cm>::}
!{(-0.8,0) }*+{\underset{a}\bullet}="1"
!{(-0.275,0)}*+{\underset{b}\bullet}="2"
!{(0.275,0)}*+{\underset{c} \bullet}="3"
!{(0.8,0) }*+{\underset{d}\bullet}="4"
"1":@/^0.7cm/"3"
"4":@/_0.7cm/"2"
 }};
(-30,-60)*{\xygraph{ !{<0cm,0cm>;<1cm,0cm>:<0cm,1cm>::}
!{(-0.8,0) }*+{\underset{a}\bullet}="1"
!{(-0.275,0)}*+{\underset{b}\bullet}="2"
!{(0.275,0)}*+{\underset{c} \bullet}="3"
!{(0.8,0) }*+{\underset{d}\bullet}="4"
"2":@/_0.7cm/"1"
"3":@/^0.7cm/"4"
 }};
(+0,-60)*{\xygraph{ !{<0cm,0cm>;<1cm,0cm>:<0cm,1cm>::}
!{(-0.8,0) }*+{\underset{a}\bullet}="1"
!{(-0.275,0)}*+{\underset{b}\bullet}="2"
!{(0.275,0)}*+{\underset{c} \bullet}="3"
!{(0.8,0) }*+{\underset{d}\bullet}="4"
"2":@/^0.5cm/"3"
"4":@/_0.7cm/"1"
 }};
(+30,-60)*{\xygraph{ !{<0cm,0cm>;<1cm,0cm>:<0cm,1cm>::}
!{(-0.8,0) }*+{\underset{a}\bullet}="1"
!{(-0.275,0)}*+{\underset{b}\bullet}="2"
!{(0.275,0)}*+{\underset{c} \bullet}="3"
!{(0.8,0) }*+{\underset{d}\bullet}="4"
"1":@/^0.7cm/"2"
"4":@/_0.7cm/"3"
 }};
(-20,-80)*{\xygraph{ !{<0cm,0cm>;<1cm,0cm>:<0cm,1cm>::}
!{(-0.8,0) }*+{\underset{a}\bullet}="1"
!{(-0.275,0)}*+{\underset{b}\bullet}="2"
!{(0.275,0)}*+{\underset{c} \bullet}="3"
!{(0.8,0) }*+{\underset{d}\bullet}="4"
"4":@/_0.7cm/"1"
"3":@/_0.5cm/"2"
 } };
(+20,-80)*{\xygraph{ !{<0cm,0cm>;<1cm,0cm>:<0cm,1cm>::}
!{(-0.8,0) }*+{\underset{a}\bullet}="1"
!{(-0.275,0)}*+{\underset{b}\bullet}="2"
!{(0.275,0)}*+{\underset{c} \bullet}="3"
!{(0.8,0) }*+{\underset{d}\bullet}="4"
"2":@/_0.7cm/"1"
"4":@/_0.7cm/"3"
 }};
(0,-100)*{\xygraph{ !{<0cm,0cm>;<1cm,0cm>:<0cm,1cm>::}
!{(-0.8,0) }*+{\underset{a}\bullet}="1"
!{(-0.275,0)}*+{\underset{b}\bullet}="2"
!{(0.275,0)}*+{\underset{c} \bullet}="3"
!{(0.8,0) }*+{\underset{d}\bullet}="4"
"3":@/_0.7cm/"1"
"4":@/_0.7cm/"2"
 }};
{\ar@{-}@/_0pc/(-4.5,-1.5);(-15.5,-16.5)}
{\ar@{-}@/_0pc/(4.5,-1.5);(15.5,-16.5)}
{\ar@{-}@/_0pc/(-20.0,-22.5);(-28.0,-31.5)}
{\ar@{-}@/_0pc/(-20.0,-22.5);(-15.0,-31.5)}
{\ar@{-}@/_0pc/(-20.0,-22.5);(+16.0,-31.5)}
{\ar@{-}@/_0pc/(+13.0,-22.5);(+23.0,-31.5)}
{\ar@{-}@/_0pc/(+11.0,-22.5);(+7.0,-31.5)}
{\ar@{-}@/_0pc/(8.0,-22.5);(-25.0,-31.5)}
{\ar@{-}@/_0pc/(-36,-42.5);(-36.0,-51.5)}
{\ar@{-}@/_0pc/(-35,-42.5);(-15.0,-51.5)}
{\ar@{-}@/_0pc/(-9,-42.5);(+12.0,-51.5)}
{\ar@{-}@/_0pc/(-11,-42.5);(-30.0,-51.5)}
{\ar@{-}@/_0pc/(+15.0,-42.5);(+15.0,-51.5)}
{\ar@{-}@/_0pc/(+12.0,-42.5);(-8.0,-51.5)}
{\ar@{-}@/_0pc/(-41,-62.5);(-36.0,-71.5)}
{\ar@{-}@/_0pc/(-40,-62.5);(-5.0,-71.5)}
{\ar@{-}@/_0pc/(-16,-62.5);(+2.0,-71.5)}
{\ar@{-}@/_0pc/(-18,-62.5);(-35.0,-71.5)}
{\ar@{-}@/_0pc/(+6.0,-62.5);(+2.0,-71.5)}
{\ar@{-}@/_0pc/(+2.0,-62.5);(-31.0,-71.5)}

{\ar@{-}@/_0pc/(-9,-82.5);(-19,-91.5)}
{\ar@{-}@/_0pc/(-37.0,-82.5);(-27.0,-91.5)}
\end{xy}$$
\end{theorem}

\remark Note that, although every minimal degeneration is of the form $W\oplus M<_{\rm mdeg}W\oplus M'$ as above, the choice of $W$ is not arbitrary, that is, addition of an arbitrary $W$ might lead to a non-minimal degeneration. The precise conditions on $W$ neccessary for this degeneration to be minimal will be described in \cite{boosphd}; as a consequence, it will be shown in \cite{boosphd} that all minimal degenerations are of codimension $1$. \\[2ex]
We have thus obtained a constructive way of describing an orbit closure $\overline{\Orb_A}$ of a $2$-nilpotent matrix $A$ in terms of its corresponding link pattern: by repeated application of the local changes to the arrows as in the theorem, we produce a list of all link patterns corresponding to matrices $A'\in N_n^{(2)}$ such that $\Orb_{A'}\subset\overline{\Orb_A}$ (although this list will contain repetitions due to non-minimal degenerations).

\section{$B_n$-orbits in arbitrary nilpotent matrices}\label{general}

In this section, we consider the action of $B_n$ on $\N_n$ by conjugation in general (for the analogous problem of the action of $B_n$ on $\mathfrak{n}_n$, see \cite{HR}).\\[1ex]
The starting point is the following observation (see \cite{Halbach}):\\[1ex]
\example Consider the action of $B_3$ on $\mathcal{N}_3$ via conjugation. Then the matrices $\left[\begin{array}{lll}0&0&0\\ 1&0&0\\ \lambda&1&0\end{array}\right]$ for $\lambda\in k$ are pairwise non-conjugate. Furthermore, on the open set $U\subset \mathcal{N}_3$ of nilpotent matrices $A=\left[\begin{array}{lll}a&b&c\\ d&e&f\\ g&h&i\end{array}\right]$ where $g\not=0$ or $dh\not=eg$, the map $U\rightarrow{\bf P}^1$ given by $A\mapsto (g:dh-eg)$ is surjective and $B_3$-invariant.\\[2ex]
We generalize some aspects of this example to arbitrary $n$.

\subsection{Generic normal form}\label{gnf}

It is appropriate to reformulate the problem as follows: we consider the action of ${\rm GL}(V)$ on pairs $(F_*,\varphi)$ consisting of a complete flag $0=F_0\subset F_1\subset\ldots\subset F_n=V$ and a nilpotent operator $\varphi\in{\rm End}(V)$ of an $n$-dimensional $k$-vector space $V$. Then the orbits of this action are precisely the orbits of $B_n$ in $\N_n$ since the variety of complete flags is isomorphic to the homogeneous space $\GL_n/B_n$.

\begin{theorem} The following properties of a pair $(F_*,\varphi)$ consisting of a complete flag and a nilpotent operator of an $n$-dimensional $k$-vector space $V$ are equivalent:
\begin{enumerate}
\item $\dim \varphi^{n-k}(F_k)=k$ for all $k=1,\ldots,n-1$,
\item $\dim(\varphi^{n-k}(F_k)+F_k)/F_k=k$ for all $k=1,\ldots,n-1$, or, equivalently, all induced maps $\varphi:F_k\rightarrow V/F_{n-k}$ are invertible,
\item there exists a unique basis $v_1,\ldots,v_n$ of $V$ such that
\begin{enumerate}
\item $F_k=\langle v_1,\ldots,v_k\rangle$ for all $k=0,\ldots,n$,
\item $\varphi(v_k)=v_{k+1}\bmod\langle v_{k+2},\ldots,v_n\rangle$ for all $k=1,\ldots,n$.
\end{enumerate}
\end{enumerate}
\end{theorem}

\begin{proof} Obviously, the second property implies the first. We show that the third property implies the second one; so assume there exists a basis $v_1,\ldots,v_n$ with the properties (a) and (b). By an easy induction, we have
$$\varphi^k(v_l)=v_{k+l}\bmod\langle v_{k+l+1},\ldots,v_n\rangle$$
for all $k+l\leq n$, and $\varphi^k(v_l)=0$ if $k+l>n$. We thus have
$$\varphi^{n-k}(F_k)=\langle\varphi^{n-k}(v_1),\ldots,\varphi^{n-k}(v_k)\rangle=\langle v_{n-k+1},\ldots,v_n\rangle,$$
and the second property follows since $F_{n-k}=\langle v_1,\ldots,v_{n-k}\rangle$.\\[1ex]
Conversely, assume that $\dim\varphi^{n-k}(F_k)=k$ for all $k$. In particular, we have $\varphi^{n-k}(V)=\varphi^{n-k}(F_k)$, and thus $\dim\varphi^{n-k}(V)=k$ and $\dim{\rm Ker}(\varphi^{n-k})=n-k$ for all $k$. We choose an arbitrary basis $w_1,\ldots,w_n$ of $V$ which is adapted to $F_*$, that is, such that $F_k=\langle w_1,\ldots,w_k\rangle$ for all $k$. Then, for all $k$, the elements $\varphi^{n-k}(w_1),\ldots,\varphi^{n-k}(w_k)$ generate the $k$-dimensional space $\varphi^{n-k}(F_k)$, thus they form a basis of this space. We can thus write the element $\varphi^{n-1}(w_1)\in\varphi^{n-k}(F_k)$ uniquely as $$\varphi^{n-1}(w_1)=\sum_{i=1}^kb_{k,i}\varphi^{n-k}(w_i),$$
and we define
$$v_k=\sum_{i=1}^kb_{k,i}w_i$$
for all $k$. Note that the elements $v_k$ do not depend on the choice of basis elements $w_1,\ldots,w_n$. We have $b_{k,k}\not=0$: otherwise $\varphi^{n-1}(w_1)=\sum_{i<k}b_{k,i}\varphi^{n-k}(w_i)$, and application of $\varphi$ yields $0=\sum_{i<k}b_{k,i}\varphi^{n-(k-1)}(w_i)$ and thus $b_{k,i}=0$ for all $i$ by linear independence of the elements $\varphi^{n-(k-1)}(w_i)$. Then $\varphi^{n-1}(w_1)=0$, a contradiction. Since the elements $w_k$ form a basis and the $b_{k,k}$ are non-zero, the elements $v_k$ form a basis, too, which is again adapted to $F_*$.\\[1ex]
We have $$\varphi^{n-k}(v_k)=\sum_{i=1}^kb_{k,i}\varphi^{n-k}(w_i)=\varphi^{n-1}(w_1)=v_n$$
by definition. For $k+l>n$, we thus have
$$\varphi^{k}(v_l)=\varphi^{k+l-n}(\varphi^{n-l}(v_l))=\varphi^{k+l-n}(\varphi^{n-1}(w_1))=0.$$
It follows that $v_{k+1},\ldots,v_n$ belong to ${\rm Ker}(\varphi^{n-k})$, thus they form a basis of this space for dimension reasons. It also follows that $\varphi^{n-k}(v_1),\ldots,\varphi^{n-k}(v_k)$ form a basis of $\varphi^{n-k}(V)$.\\[1ex]
Writing $\varphi^{k}(v_l)=\sum_{i=1}^nc_{k,l,i}v_i$, we apply $\varphi^{n-k-l}$ and calculate
$$\varphi^{n-k-l}(v_{k+l})=v_n=\varphi^{n-l}(v_l)=\sum_ic_{k,l,i}\varphi^{n-k-l}(v_i)$$
$$=c_{k,l,k+l}\varphi^{n-k-l}(v_{k+l})+\sum_{i<k+l}c_{k,l,i}\varphi^{n-k-l}(v_i),$$
and thus $c_{k,l,k+l}=1$ and $c_{k,l,i}=0$ for all $i<k+l$ by linear independence of $\varphi^{n-k-l}(v_1),\ldots,\varphi^{n-k-l}(v_{k+l})$. We thus have
$$\varphi^k(v_l)=v_{k+l}+\sum_{i>k+l}c_{k,l,i}v_i$$
for all $k+l\leq n$, and, in particular,
$$\varphi(v_k)=v_{k+1}+\sum_{i>k+1}c_{1,k,i}v_i$$
for all $k$. The basis $v_1,\ldots,v_n$ thus has the claimed properties.\end{proof}

For $0\leq a,b\leq n$ and a matrix $A\in \N_n$, define $A_{(a,b)}$ as the submatrix formed by the last $a$ rows and the first $b$ columns of $A$.

\begin{corollary}\label{hnf} The following conditions on a matrix $A\in\N_n$ are equivalent:
\begin{enumerate}
\item for $k=1,\ldots,n-1$, the first $k$ columns of $A^{n-k}$ are linearly independent,
\item for $k=1,\ldots,n-1$, the minor $\det ((A^{n-k})_{(k,k)})$ is non-zero,
\item $A$ is $B_n$-conjugate to a unique matrix $H$ such that $H_{i,j}=0$ for $i\leq j$ and $H_{i+1,i}=1$ for all $i=1,\ldots,n-1$.
\end{enumerate}
\end{corollary}

\begin{proof} We apply the previous theorem to the vector space $V={\bf k}^n$ with coordinate basis $e_1,\ldots,e_n$, the standard flag defined by $F_k=\langle e_1,\ldots,e_k\rangle$ and the endomorphism $\varphi$ given by multiplication by $A$. The first property of the theorem immediately translates into linear independence of column vectors, whereas the second property translates to the non-vanishing of minors. The basis $v_1,\ldots,v_n$ of the theorem yields an upper-triangular base change matrix, and representing $A$ with respect to this basis yields the desired $B_n$-conjugate $H$.\\[1ex]
\end{proof}

The conditions of Corollary \ref{hnf} define an open subset of $\N_n$; we have thus found a generic normal form for nilpotent matrices up to $B_n$-conjugacy.

\subsection{Semiinvariants}\label{si}

We construct a class of determinantal $B_n$-semiinvariants on $\N_n$, that is, regular functions $D$ on $\N_n$ such that $D(gAg^{-1})=\chi(g)D(A)$ for all $g\in B_n$ and $A\in \N_n$; here $\chi$ is a character on $B_n$ called the weight of $D$. For $i=1,\ldots,n$, we denote by $\omega_i:B_n\rightarrow{\bf G}_m$ the character defined by $\omega_i(g)=g_{i,i}$; the $\omega_i$ form a basis for the group of characters of $B_n$.\\[1ex]
Fix non-negative integers $a_1,\ldots,a_s,b_1,\ldots,b_t$ such that $a_1+\ldots+a_s=b_1+\ldots+b_t=:k\leq n$. Moreover, fix polynomials $P_{i,j}(x)\in k[x]$ for $i=1,\ldots,s$ and $j=1,\ldots,t$, and denote the datum $((a_i)_i,(b_j)_j,(P_{i,j})_{i,j})$ by $P$. For all such $i$ and $j$, consider the $a_i\times b_j$-submatrices $P_{i,j}(A)_{(a_i,b_j)}$ as defined in the previous section, and form the block matrix $A^P=(P_{i,j}(A)_{(a_i,b_j)})_{i,j}$; this is a $k\times k$-matrix.

\begin{proposition} For every datum $P$ as above, the function associating to a matrix $A\in\N_n$ the determinant $\det(A^P)$ defines a $B_n$-semiinvariant regular function $D^P$ of weight $\sum_i(\omega_{a_i}+\ldots+\omega_n)-\sum_j(\omega_1+\ldots+\omega_{b_j})$ on $\N_n$.
\end{proposition}

\begin{proof} For $g\in B_n$ and $1\leq a,b\leq n$, denote by $g_{(\geq a)}\in B_a$ (resp. by $g_{(\leq b)}\in B_b$) the submatrix formed by the last $a$ rows and columns (resp. by the first $b$ rows and columns) of $g$. With these definitions, it follows immediately that
$$(gAg^{-1})_{(a,b)}=g_{(\geq a)}A_{(a,b)}g_{(\leq b)}^{-1}.$$
This yields the following equalities of block matrices
$$(gAg^{-1})^P=(P_{i,j}(gAg^{-1})_{(a_i,b_j)})_{i,j}=((gP_{i,j}(A)g^{-1})_{(a_i,b_j)})_{i,j}$$
$$=(g_{(\geq a_i)}P_{i,j}(A)_{(a_i,b_j)}{g_{(\leq b_j)}}^{-1})_{i,j}=(\delta_{i,j}g_{(\geq a_i)})_{i,j} A^P (\delta_{i,j}{g_{(\leq b_j)}}^{-1})_{i,j},$$
and thus
$$D^P(gAg^{-1})=\det((gAg^{-1})^P)=\prod_i\det(g_{(\geq a_i)})\prod_j\det(g_{(\leq b_j)})^{-1}D^P(A).$$
\end{proof}

With the aid of these semiinvariants, we can see that the entries of the normal form $H$ associated to a matrix $A$ fulfilling the conditions of Corollary \ref{hnf} depend polynomially on $A$, by describing them as the value of a special semiinvariant $D^P$:

\begin{lemma} For $i$ and $j$ such that $1\leq j\leq n-2$ and $j+2\leq i\leq n$, consider the datum $P$ as above defined by $a_1=j-1$, $a_2=i$, $b_1=j$, $b_2=i-1$, $P_{1,1}(x)=x^{n-j+1}$, $P_{1,2}(x)=0$, $P_{2,1}(x)=x$, $P_{2,2}(x)=x^{n-i+1}$. Then, for a matrix $H$ in the form of Corollary \ref{hnf}, we have $D^P(H)=H_{i,j}$.
\end{lemma}

\begin{proof} By a direct calculation, the matrix $H^P$ consists of the blocks
$$(H^P)_{1,1}=\left(\begin{array}{cccc}1&&0&0\\ &\cdots&&\vdots\\ *&&1&0\end{array}\right),\, (H^P)_{1,2}=0,$$ $$(H^P)_{2,1}=H_{(i,j)},\, (H^P)_{2,2}=\left(\begin{array}{ccc}0&\ldots&0\\ 1&&0\\ &\vdots&\\ *&&1\end{array}\right).$$
Thus, the matrix $H^P$ is lower triangular, all diagonal entries being $1$ except the $(j,j)$-entry, which equals $H_{i,j}$.
\end{proof}

It seems likely that the semiinvariants $D^P$ generate the ring of all semiinvariants at least for a certain cone of weights. The generic normal form of Corollary \ref{hnf} allows to find identities between the $D^P$ by evaluation on matrices $H$ in normal form.

\end{document}